\documentclass[final,5p,times,twocolumn]{elsarticle}

\usepackage{amssymb}
\usepackage{graphicx,amsmath}
\usepackage{ifthen}
\usepackage{endnotes}
\def\beq{\begin{equation}}
\def\eeq{\end{equation}}
\def\bea{\begin{eqnarray}}
\def\eea{\end{eqnarray}}
\def\nn{\nonumber}

\biboptions{comma,round,authoryear}

\journal{Energy Policy}

\begin{document}

\begin{frontmatter}


\title{FTT:Power : A global model of the power sector with induced technological change and natural resource depletion}
\author{Jean-Fran\c{c}ois Mercure \corref{cor1}}
\ead{jm801@cam.ac.uk}
\cortext[cor1]{Corresponding author: Jean-Fran\c{c}ois Mercure}
\address{Cambridge Centre for Climate Change Mitigation Research (4CMR), Department of Land Economy, University of Cambridge, 19 Silver Street, Cambridge, CB3 1EP, United Kingdom}

\begin{abstract}
This work introduces a model of Future Technology Transformations for the power sector (FTT:Power), a representation of global power systems based on market competition, induced technological change (ITC) and natural resource use and depletion. It is the first component of a family of sectoral bottom-up models of technology, designed for integration into the global macroeconometric model E3MG. ITC occurs as a result of technological learning produced by cumulative investment and leads to highly nonlinear, irreversible and path dependent technological transitions. The model uses a dynamic coupled set of logistic differential equations. As opposed to traditional bottom-up energy models based on systems optimisation, such differential equations offer an appropriate treatment of the times and structure of change involved in sectoral technology transformations, as well as a much reduced computational load. Resource use and depletion are represented by local cost-supply curves, which give rise to different regional energy landscapes. The model is explored for a single global region using two simple scenarios, a baseline and a mitigation case where the price of carbon is gradually increased. While a constant price of carbon leads to a stagnant system, mitigation produces successive technology transitions leading towards the gradual decarbonisation of the global power sector. 

\end{abstract}


\end{frontmatter}


\section{Introduction}

The future level of anthropogenic greenhouse gas (GHG) emissions is one of the primary unknowns in estimating the rate of climate change in the medium and long term. Different assumptions about emissions pathways lead to widely differing warming temperatures, ranging between 2$^\circ$C, for high mitigation policies to 6-8$^\circ$C for the most pessimistic non mitigation scenarios \citep{IPCCAR4}. Anthropogenic GHG emissions depend on human activities, and thus on the structure of the future economic system. They stem primarily from energy use (57\%) and land use change (17\%) \citep[figures for 2004]{IPCCAR4}. While the second is a complex subject to model (see for instance \cite{IMAGE}), the first is better understood and different approaches have been used to model  energy demand and supply \citep{IEAWEO2010, IMAGE, MESSAGE, MARKAL}. The problem of GHG emissions mitigation requires changes to be made to the structure of global energy use. Since the latter lies at the very core of the world's economy, these changes have deep implications and effects felt in all aspects of society. Therefore, a simulation of GHG emissions cannot easily be separated from simulations of the global economic system. One approach to explore Energy-Economy-Environment (E3) interactions was made with an input-output structure based macroeconometric model of the global economy E3MG \citep{Barker2006, Dagoumas2010, Scrieciu2010}. E3MG is a disaggregated model that features 20 world regions, 12 energy carriers, 19 energy users, 28 energy technologies, 14 atmospheric emissions and 42 industrial sectors.

The decarbonisation of the global power system depends first and foremost on the rate at which highly emitting technologies based on fossil fuels can be substituted for cleaner ones. While fossil fueled electricity generation technologies are mature and well established, cleaner systems such as renewable energies are more expensive and therefore not yet competitive. Policy favouring green technologies are expected to play an important role in the transition towards low GHG emissions in the energy sector. The costs of technologies tend to decrease with cumulative investment, an effect termed technological learning, and moreover, newer technologies see their costs decrease more rapidly than those of mature technologies \citep{IEAlearning, Pan2007, McDonald}. R \& D investments and government subsidies to green technologies are expected to be necessary only for an initial `push', until their reducing costs reach competitive levels, giving rise to technological change \citep{IEAlearning}. The composition of future energy landscapes are highly dependent on the development in time of policy portfolios, where some strategies can generate technological avalanches which would not have occurred otherwise. Energy landscapes are, however, very local in nature and technological substitution occurs according to what interactions are possible in a particular world region. The energy portfolio of each region is in large parts determined by the availability of natural resources. Therefore, any modelling attempt of technological substitution should reflect the local nature of natural resources.

This paper introduces the model FTT:Power (Future Technology Transformations in the Power sector), a dynamic model of the global power sector based on market price competition, technological substitution and resource use and depletion. It is the first member of a family of bottom-up dynamic technology substitution models, built as components for E3MG. FTT:Power makes use of a dynamic set of coupled logistic differential equations, similar to those used to represent population growth \citep{Verhulst1838} and competition in biological systems \citep{Lotka1925,Volterra1939}. Evolutionary dynamics have been used to describe technological change theoretically \cite[see for instance][]{Dercole2008}. Technology transitions have been shown numerous times to follow empirically logistic functions \citep[for a review see][]{Grubler1999}, in particular in energy systems \citep{Marchetti1978}. There exists an extensive literature where logistic systems of equations have been used in the analysis of growth and competition in markets \citep[for instance][]{Bass1969,Sharif1976,Bhargava1989,Morris2003}. While these systems follow closely logistic functions, their underlying changeover timescales have never been explained satisfactorily by any theory. The subject has been explored more recently by \cite{Anderson2007}, where they introduced a family of coupled differential logistic equations to represent energy technology transitions. Their system is not dynamic, however, and leads to static results when technology costs do not change, and furthermore lack appropriate symmetry properties. Meanwhile, current energy technology substitution models generally fail to reproduce the sigmoid ($S$-shaped) character of technological transitions that stem from their dynamic nature. Most models use instead procedures of energy systems optimisation based on cost minimisation, which produces optimal equilibrium solutions. They lack the intrinsic time constants (of the order of 50 to 100 years) associated with technological transitions, and smoothing of changes in time are generally imposed exogenously. 

The FTT:Power model is the first power sector model based on a dynamic set of coupled logistic differential equations, of the Lotka-Volterra family, which represents gradual technological substitution processes. It is parametrised by natural time constants which determine the rates at which technological transitions can occur. It uses a probabilistic treatment of investor decision-making using local distributions of cost values and resulting likelihood of technological switching. It allows for non-rational investor behaviour, assuming that the trend rather than individual actions is oriented towards choosing lowest cost technologies. Due to the long time constants involved, which are associated with lifetimes and construction times, however, the model does not reach an equilibrium state as it would in an optimisation of the system for cost minimisation, unless all parameters were to remain stable for periods several times longer than these time constants. Due to the simplicity of evaluating in time steps a set of differential equations, FTT:Power takes a small fraction of the computational time required by optimisation models; it is thus likely to be one of the most compact and fastest models of technology around. Additionally, while technology learning curves generate multiple optimal points to optimisation models and associated computational difficulties, they do not require additional computational time in this model. 

The model is designed to calculate in parallel the dynamic evolution of local energy landscapes in 20 world regions based on local cost values. These are constructed using a standard framework which includes components representing policy decisions such as taxes and carbon markets, experience curves and the availability of natural resources and their cost. The treatment of natural resource use and depletion is inspired by that featured in the TIMER model, which defines global sets of local cost-supply curves \citep{IMAGE, HoogwijkThesis, Hoogwijk2004, Hoogwijk2005, Hoogwijk2009, Rogner1997}, but expanded to include a representation of uncertainty and evolution in resource assessments. As an integral part of E3MG, FTT:Power is designed to reproduce feedback interactions between the global economy and the energy sector, which occurs through the energy demand, itself highly influenced by energy prices. Electricity prices are derived through FTT:Power with the marginal costs of each technology, which depend on those of natural resources and technology. 

This paper is divided into three sections. The logistic model of technology substitution in differential equation form and the associated probabilistic approach to local cost distributions are first introduced. We discuss the effects of technological learning onto cost values, and thus to decision making, given local policy portfolios. We then describe how this equation is constrained using a simple representation of technological limits imposed by local grid networks. The second section introduces our framework of natural resource use and depletion featuring an explicit representation of uncertainty in global resource distributions and availability. This is finally followed by the presentation of one set of global results in order to explore and understand the properties of FTT:Power. For this purpose, the model was run for this paper by itself without feedback with E3MG, in order to avoid the significant increase in complexity of the model when feedback effects are present. The results were obtained using two sets of simplified global policy portfolios, a baseline and a mitigation scenario. We discuss the appearance of a phenomenon we name the energy technology ladder, which emerges from the complexity of the set of equations at the root of FTT:Power. This effect is a clear representation of irreversible and path dependent induced technological change and the way in which the energy sector might evolve in high mitigation scenarios in order to achieve high levels of decarbonisation efficiently at low costs. This result stems from our dynamic approach rather than being imposed through a cost minimisation procedure and associated assumption of market equilibrium, and is consistent with the non-equilibrium Post-Keynesian approach of E3MG. Additional information regarding FTT:Power can be found in \cite{Mercure2011TWP}.

\section{Dynamics as a set of differential equations \label{sect:dynamics}}

\subsection{The shares equation}

\begin{figure}[hbt]
	\begin{center}
		\includegraphics[width=1\columnwidth]{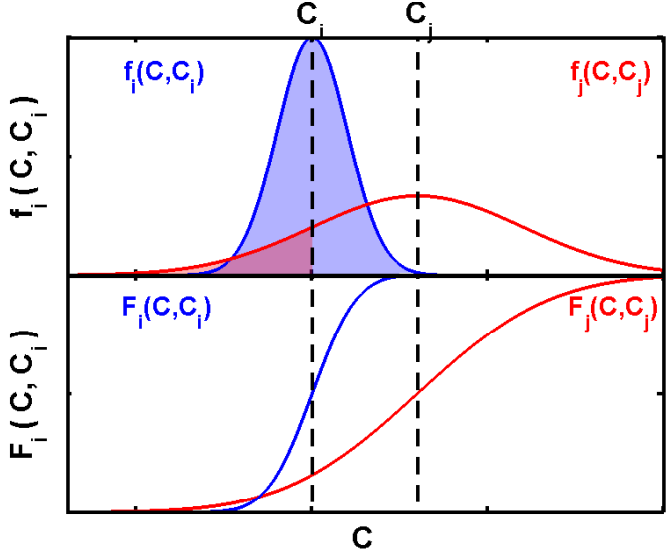}
	\end{center}
	\caption{$Top$ Probability distributions $f_i(C,C_i)$ and $f_j(C,C_j)$ of cost values for two technologies $i$ and $j$ based on real recently experienced costs. The number of units of technology $j$ that come out at a cost cheaper than the median value of the distribution of technology $i$ corresponds to red shaded area, a value much smaller than that of the reverse, which is the blue area. $Bot.$ Cumulative probability distribution functions for $i$ and $j$.}
	\label{fig:Figure1}
\end{figure}

We introduce the central assumption of our model by presenting a probabilistic framework designed to represent the results of local decision-making by investors in the electricity sector. The electricity sector sees a local varying demand for power to which it reacts by attitudes regarding the building of new or maintaining existing power stations. On the ground, investors are faced with a choice between several avenues, and they perform a comparison of the Levelised Cost of Electricity (LCOE, defined in section \ref{sect:LCOE}). This cost, however, is not strictly well defined, since a large number of aspects of very local nature affect the cost of specific projects, such as the length of power lines or the cost of land, and decisional issues which may or may not lead to a rational choice of the option with the lowest LCOE. We therefore take costs as being distributed when examined from an aggregated viewpoint at the regional or national level, in order to avoid looking into the details of these local aspects. Therefore, even if one particular technology is on average less expensive than another, there are almost always individual cases where it is the reverse that is true, and therefore a probabilistic treatment is appropriate. Thus we assume that it is the trend rather than individual decisions which leans towards technologies with lower LCOEs.

We define $f_i(C, C_i)$ as a probability distribution of the cost of technology $i$ based on real recent instances of the cost of new units (taken from \cite{IEAProjCosts}), which has a median value $C_i$. The probability, when building a new unit, that it comes at a cost lower than the value $C$ is the cumulative probability distribution $F_i(C,C_i)$. For two technologies $i$ and $j$, the fraction of new units of $j$ whose cost turn out lower than the median value of $i$ is $F_j(C_i,C_j)$, while the fraction of units of $i$ whose cost is lower than the median of $j$ is $F_i(C_j,C_i)$. These can also be expressed in terms of the cost difference $\Delta C_{ij}$. This is shown in figure~\ref{fig:Figure1}. Using these arguments, the differential shares equation can be constructed.

We model technology switching using a pairwise comparison of costs, and flows of units of market shares $S_i$ of different technologies, that is, the shares of electricity generation capacity. We will demonstrate that this is equivalent to a simultaneous comparison of the costs of all options. We first determine the likely number of energy production units flowing from category $j$ towards category $i$, denoted $\Delta S_{j \rightarrow i}$. The number of units of market share that can flow away from technology $j$ is a fraction of the existing fleet, and therefore is proportional to the market share itself $S_j$. The rate, therefore, at which power plants of category $j$ will be closed down, if none are scrapped before the end of their lifetime $\tau_j$, is $S_j / \tau_j$. Meanwhile, the rate at which technology $i$ can grow is also proportional to its share of the market. This is similar to population growth: the technology is assumed to `breed' itself, an assumption which reflects the growing ability of an industry to expand as its sales increase, and forms the basis of the logistic based market competition literature \citep{Bass1969,Sharif1976,Bhargava1989,Morris2003}. The rate at which technology $i$ can grow is inversely proportional to the technology specific building capacity expansion time constant $t_i$, and is therefore $S_i / t_i$. In a unit of time $dt$, where $t$ is measured in years throughout, we assume that the number of units flowing from $j$ to $i$ is proportional to the probability that the cost of $i$ is less than that of $j$, $F_i(C_j,C_i)$, and conversely, leading to
\bea
\Delta S_{j \rightarrow i} &\propto& {S_j \over \tau_j}{S_i \over t_i} F_i(\Delta C_{ij}) \Delta t\nn\\
\Delta S_{i \rightarrow j} &\propto& {S_i \over \tau_i}{S_j \over t_j} F_j(\Delta C_{ji}) \Delta t.\nn
\eea
These aggregate movements are completely independent of each other: they stem from every individual substitution event where investors make a specific choice and where costs differ for local reasons which we ignore here. There results a net flow, which is the difference between these two terms, denoted $\Delta S_{ij}$, and represents every substitution event that occurs during the time interval $\Delta t$. In order to find the total number of units flowing towards $i$ from all other technologies, we perform the sum of $\Delta S_{ij}$ over $j$, which yields the shares equation central to this model:
\beq
\Delta S_i = \sum_j S_i S_j \left( A_{ij} F_i(\Delta C_{ij})-A_{ji} F_j(-\Delta C_{ij})\right)\Delta t,
\label{eq:Shares}
\eeq
where $A_{ij} = K/\tau_i t_j$ is a matrix of substitution frequencies, which is not symmetrical, with $K$ a time scaling constant.\endnote{Note that for cases where substitution occurs faster than the lifetime of a plant, exceptions can be inserted in $A_{ij}$. This may occur in the case of, for instance, the retrofitting of carbon capture and storage to an existing plant, or fuel switching from coal to biomass.} This matrix contains all the natural time constants associated with induced technological change. Its detailed derivation involves a lengthy calculation that cannot be included here but will be the subject of a forthcoming publication. This theory yields a changeover (transition) time, which depends on the strength of investor preferences $F_i(\Delta C_{ij})$, and ranges between the shortest possible decommission time $\tau$ ($F_i(\Delta C_{ij}) = 1$) and infinity ($F_i(\Delta C_{ij}) = 0$). Note that $A_{ij}$ is $not$ symmetrical, and this stems from the fact that rates of technology uptake and rates of decommission are not identical across technologies. For example, technology $i$ being decommissioned at rate $\tau_i^{-1}$ and replaced by technology $j$ at a rate proportional to $t_j^{-1}$ does not occur at the same rate as the reverse process where technology $j$ is decommissioned, at rate $\tau_j^{-1}$, and replaced by $i$ at rate proportional to $t_i^{-1}$.

Additionally, it can be seen that the equation is invariant over permutations of $i$ and $j$, required in order to have $\sum_{ij} dS_{ij} = 0$ and $\sum_i S_i = 1$. It can be shown numerically that the result does not depend on the order over which the calculation is performed. This means that, if the same result is obtained irrespective of the order over which the pairwise comparisons are done, that the total operation corresponds to a simultaneous comparison of all technologies. This is ensured by the differential form of eq. (\ref{eq:Shares}), where exchanges are performed over infinitesimal time intervals, and is further demonstrated with examples in section \ref{sect:examples}.

\subsection{Induced Technological Change}
One of the most important drivers of change in the electricity sector is the cost reduction that stems from technological learning, and have been included in several energy models, notably in the preceding energy submodel of E3MG \citep{Kohler2006b}. It is well established that the repetitive production of goods gives rise to improvements in production methods and economies of scale that lead to cost reductions, and these depend directly on the number of units produced~\citep{Arrow1962, Grubler1999}. Moreover, the functional form of the experience curve is very simple, it is that of a power law of which the negative exponent, denoted $-b_i$ here, is related to the progress rate \citep{IEAlearning, Berglund2006}. Thus, the logarithm of the cost of production $C_i$ relative to an initial value depends linearly on the logarithm of the total number of units $W_i$ sold since the very first one came out of the factory. In linear terms, the experience curve is expressed as 
\beq
C_i(t) = C_{0,i} \left({W_i(t) \over W_{0,i}}\right)^{-b_i},
\label{eq:ExpCurve}
\eeq
where $C_{0,i}$ is the cost associated with the cumulative number of units $W_{0,i}$ produced up to an arbitrary starting time $t=0$. Values for $W_0$ are not trivial to find in the cases of old technologies since they involve the number of units currently operating but also the total number of units which have been abandoned or demolished; learning is however minimal in such cases. Where learning is important, with newer technologies, these correspond approximately to the current number of  installed units since most of them are still in operation. The learning process can also be seen in an incremental way, through the accumulation of knowledge, expressed as
\beq
\Delta C_i = -b_i {C_i \over W_i} \Delta W_i,
\eeq
which involves implicitly starting values $C_{0,i}$ and $W_{0,i}$, and is equivalent to the previous equation.
 
Technology categories for learning and those represented in a model may not necessarily coincide, and thus a certain amount of mixing, or knowledge $spillover$, may have to be included. In other words, particular sets of categories may be closely related technologically and a learning spillover matrix $B_{ij}$ should be defined in order to calculate an effective $W_i$ from incremental positive capacity additions:
\beq
W_i(t) = \sum_j B_{ij} \left\{ \begin{array}{ll} \int_0^t\left( {dU_j(\tau) \over d\tau} + \delta_j U_j(\tau) \right)d\tau,  & {dU_j(\tau) \over d\tau} > 0\\ 
\int_0^t  \delta_j U_j(\tau)d\tau, & {dU_j(\tau) \over d\tau} \leq 0 \end{array}, \nn\right.
\eeq
where $U_j$ is the capacity of technology $j$ (defined below in eq.~\ref{eq:U}). This equation thus insures that knowledge is shared between related technologies.\endnote{Technologies with knowledge spillover include for instance coal and biomass gasification, offshore and onshore wind, combined cycle gas turbines (CCGT) and integrated gasification combined cycle (IGCC). These connexions can arise for instance through the use of similar mechanical parts that involve similar production methods, susceptible to economies of scale.} There exists an extensive literature on learning and progress rates for all sorts of goods beyond the electricity sector (see for instance \cite{Kohler2006, Pan2007, Grubler1999}). Within the power sector, learning rates have been compiled in both \cite{IEAlearning} and \cite{McDonald}.

The consequences of including experience curves are very important but lead to strong non-linearities in any model due to their strong associated concavity \citep{Berglund2006}, and to systems with multiple optimum values or none altogether \citep{Grubb2002}, making optimisation models difficult to solve, a problem that does not arise here. They moreover produce a strong path dependence, an effect we consider fundamental to the process of technological change. Since cost reductions depend on the ratio of $W_i/W_{0,i}$ according to an inverse power law, new technologies see higher cost changes than mature ones for identical numbers of additional sales. New expensive technologies have the potential to become less expensive than mature ones given that enough units are produced. Therefore, in a market where investors choose according to costs, new technologies require government R\&D investment and subsidies in order to become competitive, so-called market push policies \citep{Nemet2009}. However, with learning, such technologies may in time pass a competitiveness threshold, after which an avalanche effect can happen, where learning results in more units sold, which results in more learning and so on, and produces an induced technology transition. No single optimum solution exists since costs depend on paths, and an unlimited number of paths can be followed, which depend primarily on the assortment of possible policy decisions. These effects are inherently dynamic in nature, and are thus appropriately represented by a dynamic set of logistic differential equations connected to experience curves, which lead naturally to a behaviour very close to the observed $S$-shaped technological transitions, as explored for instance by \cite{Grubler1999}.

\subsection{Levelised cost of electricity \label{sect:LCOE}}
Investors in the power sector face complex decisions involving a very large number of parameters. However, a standard framework exists which is used by most of the industry worldwide, that of the LCOE. It is not the focus in this work to review this framework, and we refer the reader to the recent report by the \cite{IEAProjCosts} on the costs of various electricity generation technologies. We thus use the LCOE in the following form,
\beq
LCOE_i(t) = {\sum_{t = 0}^{\tau_i} {TI_i(t) + OM_i(t) + FC_i(t) + CC_i(t) \over (1+r)^t} \over \sum_{t = 0}^{\tau_i} {EP_i(t) \over (1+r)^t}},
\label{eq:LCOE}
\eeq
where $TI_i$ denotes the specific technology investment cost, $OM_i$ the operation and maintenance costs, $FC_i$ the fuel costs, $CC_i$ the carbon cost component associated with emissions allowances or taxes where applicable, $r$ the technology or region specific discount rate and $EP_i$ is the energy that the power station is expected to produce. The LCOE is the cost which is compared in the shares equation. Since all components of the LCOE can be taken as distributions, the LCOE is also itself a distribution.

\subsection{Complete set of variables}
The evolution of the shares $S_i$, as given above, correspond to changes in the composition of the electricity generation $G_i$, in GWh,
\bea
\label{eq:G}
G_i(t) &=& U_i(t) CF_i(t) = U_{tot} S_i(t) CF_i(t),\\
U_{tot} &=& \sum_i U_i,\nn
\eea
where, $U_i$ is the capacity in units of GW, $CF_i$ is the capacity (or load) factor, defined as the average fraction of time that a unit runs at full output, times a constant of conversion between $G_i$ and $U_i$ of 8760 h/y. The capacity $U_i$ can therefore be expressed as
\bea
\label{eq:U}
U_i(t) &=& {S_i(t) D(t) \over \overline{CF}(t)},\\
\overline{CF}(t) &=& \sum_i S_i(t)CF_i(t),\nn
\eea
where $D(t)$ is the electricity demand at time $t$ and $\overline{CF}$ is the weighed average of the capacity factor over the whole electricity sector. Changes in capacity stem from independent changes in three variables, $S_i$, $D$ and $\overline{CF}$, emphasising three possible processes. This can be expressed in differential form,
\beq
dU_i = {S_i \over \overline{CF}} dD + {D \over \overline{CF}} dS_i - {S_i D \over \overline{CF}^2}d\overline{CF},
\eeq
where the first term expresses changes in capacity due to changes in electricity demand, the second due to changes in the composition of the electricity sector, and the third due to changes in the efficiency at which the electricity sector is used.

Investment $I(t)$ in new generation capacity correspond to positive changes in $U_i$ times their cost $C_i$:
\beq
I_i(t) = \left\{ \begin{array}{ll}C_i(t)\left( {dU_i(t) \over dt} + \delta_i U_i(t) \right),  & {dU_i(t) \over dt} > 0\\ 
C_i(t) \delta_i U_i(t), & {dU_i(t) \over dt} < 0\end{array} \right.
\eeq
where $\delta_i$ is the rate of power plant decommissioning equal to the inverse of the lifetime $\tau_i$.

Finally, from $G_i(t)$ and emissions factors obtained from the IPCC Guidelines \citep{IPCCGuidelines}, it is straightforward to calculate the GHG emissions that stem from electricity production. Defining $\alpha_i$ as the emission factors, GHG emissions and cumulative emissions are
\bea 
E(t) &=& \sum_i \alpha_i G_i(t),\\
E_{tot}(t) &=& \int_0^t \sum_i \alpha_i G_i(t') dt'.
\eea
The carbon cost component of the LCOE is proportional to $\alpha_i$ times the price of carbon, where we neglect efficiency improvements within technology categories.

\subsection{Technical Constraints}
\begin{figure}[hbt]
	\begin{center}
		\includegraphics[width=1\columnwidth]{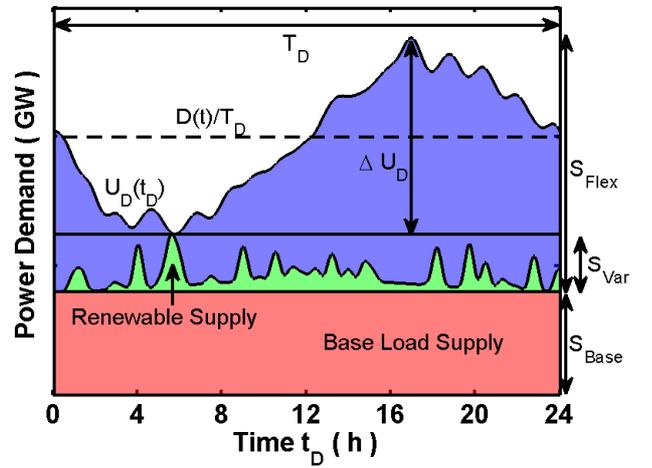}
	\end{center}
	\caption{Sketch of a hypothetical profile of daily power demand $U_D(t_D)$ as a function of the time of day $t_D$, expressed in capacity units (GW), and how it might be met by various types of supply. Areas represent amounts of energy. The green area is the supply of variable renewable energy, while the red rectangle is the base load supply. The difference between the sum of these and the demand must be met by flexible systems that can vary their output from zero to a large fraction of the total demand, shown as a blue area. $D(t)/T_D$ is the average daily power demand and $T_D$ is the length of a day. Double arrows on the right correspond to the shares of capacity.}
	\label{fig:Figure2}
\end{figure}
The shares equation treats every technology equally. However, real world electricity grids face complex optimisation problems due to fluctuations in time of both demand and supply. While the demand varies according to consumer habits, the supply fluctuates according to the varying nature of some natural resources such as wind and solar irradiation. For grid stability, demand and supply must be met at every second, and therefore some flexibility must exist to make both meet. This is done with power plants which have the ability to vary their output rapidly, such as combined cycle gas turbines (CCGT), diesel generators or hydroelectric dams. However, a large fraction of the electricity generation comes from power stations which do not have the ability to vary their output, but have constant power production (i.e. base load systems), such as nuclear reactors or coal power stations, or have an uncontrollable variable output, such as solar panels or wind turbines. Figure~\ref{fig:Figure2} sketches a typical daily demand and supply situation, where in green is represented the variable renewables contribution, while in red is shown the base load generation. $U_D$ is the power demand during a day, and therefore energy units correspond to areas. We observe that the area in blue, which during a day may vary between zero and a large fraction of the total output, must be covered by flexible energy systems with fast slew rates. These flexible systems can, however, also contribute to base load production.

We define here very simple rules that appropriately captures the relevant dynamics without requiring complex hourly based simulations of the power system merit order. We first classify generation technology into three categories, $S_{Base}$, $S_{Var}$ and $S_{Flex}$, which correspond in order to the base load technologies, the variable renewables and the flexible systems. In real national power systems, the hierarchical choice of power stations follows the so-called merit order principle (see for instance \cite{IEAProjCosts}), where power units are invoked in order of cost as the demand grows, starting with base load power stations, costly diesel units being reserved for relatively rare demand spikes. We require here, however, to avoid simulating such a system and replace it with a set of simple rules. Flexible systems are necessary in order to maintain grid stability, but may in some cases incur higher operating costs than base load systems. They are in general however rewarded by higher electricity prices. The total share of flexible systems can be constrained to remain above a certain necessary value that depends on the amount of available storage, the size of the peak load demand and the amount of renewables present in the power system. Within this necessary share, a choice of flexible technologies exist between which exchanges can occur, within a restricted electricity market that meets peak load demand. Similarly, the total share of variable renewable systems can be made to remain below a certain limit that is given by the amount of available storage, the size of the peak load demand and the total capacity of flexible systems. Changes in the capacity of flexible systems can allow expansions of renewables and, conversely, reductions in the capacity of variable renewables enable reductions in the capacity of flexible systems.

In order to generate these properties, we define as $\Delta U_D$ the height of the daily peak demand function $U_D(t_D)$ averaged yearly, where $t_D$ is the time of day, relative to the lowest daily demand value. We assume that the demand $D(t)/T_D$ is the average daily power demand at year $t$. Therefore, if enough storage was available to cover any area above $D(t)/T_D$, denoted $\Delta D$, no flexible systems would be needed, but it is unlikely in most regions of the world given the large scale of this amount of energy.\endnote{As an example, using a rough calculation for the UK only based on data by \cite{Poyry2009}, levelling out the demand would require the daily storage of about 50 GWh of energy with production capacity of 8 GW, values much above the total national hydroelectric production capability.} Note that changes in available storage, distributed generation, demand response effects and smart grids, which can influence the issues treated in this section through efficiency savings, are considered outside of the scope of the present work, but their effect can be included exogenously or dynamically from another submodel through a changes in the demand $D$ or in the parameters $\Delta U_D$, $\Delta D$, $U_s$ and $E_s$. 

We thus assign to the storage generation capacity the variable $U_s$ (in GW) and to the energy storage the variable $E_s$ (in GWh). The constraint is then expressed as three inequalities:
\beq
S_{Flex}CF_{Flex} + S_{Var}CF_{Var} \geq \overline{CF} \left( {\Delta D \over D} + {U_{Var} T_D \over D}- {E_s \over D}\right),\nn
\eeq
\beq
S_{Flex} - S_{Var} \geq \left({\Delta U_D \over U_{tot}} - {U_s \over U_{tot}}\right),\nn
\eeq
\beq
S_{Base} + S_{Var} \leq \left(\overline{CF} - {1 \over 2}{\Delta U_D \over U_{tot}} + {U_s \over U_{tot}}\right),
\label{eq:Ineqs}
\eeq
where $D$ is the yearly average demand, the same value as in eq.~\ref{eq:U}, and $CF_{Flex}$ (and similarly for $CF_{Var}$) is calculated as an average of the capacity factor over all flexible systems:
\beq
CF_{Flex} = {1\over S_{Flex}}\sum_{i = Flex} S_i CF_i .\nn
\eeq
The first these inequalities expresses the constraint on amounts of energy, where flexible systems are required to produce the blue area in figure~\ref{fig:Figure2}, while the second expresses the requirement in terms of minimum capacity of flexible systems, which must cover at least $\Delta U_D$ plus the complete capacity of variable systems for when the natural resource is not available (i.e. for instance when the wind does not blow). The requirement on shares of flexible capacity given by the second inequality is indicated with the double arrow on the right of figure~\ref{fig:Figure2}. The third inequality refers to the maximum allowed capacity for base load and variable systems, which cannot exceed the area below the lowest daily demand value.\endnote{For symmetry, a fourth, trivial, constraint exist, which requires energy generation by the three types of technologies to equal total daily demand, already met by the assumptions.}

These inequalities can be transformed into limits for the value of shares. They can be interpreted in several ways. For instance, we may ask what is the minimum share of capacity allowed of flexible systems $S_{Flex}$ given the amount of storage $U_s$ and the share of capacity of variable systems $S_{Var}$. But we may also ask what is the maximum for $S_{Var}$ given $S_{Flex}$ and the other parameters. This is summarised by the following definitions for the limits on share values (denoted with a hat):
\beq
\hat{S}_i = \pm\left[\left({\Delta U_D \over U_{tot}} - {U_s \over U_{tot}}\right) + S_{Var} - S_{Flex}\right] + S_i,\nn
\eeq
\beq
\hat{S}_i = \left[\left(\overline{CF} - {1 \over 2}{\Delta U_D \over U_{tot}} + {U_s \over U_{tot}}\right) - S_{Base} - S_{Var}\right] + S_i,
\label{eq:Sharelimits}
\eeq
where in the first equation the first term is positive if $i$ refers to flexible output and negative if it refers to variable renewables \endnote{The last term $Si$ is present in order to take out the $i$th contribution, i.e. the contribution of the technology in question.}. The most constraining of both limits is taken, which can change according to the situation. The constraints summarised by eq. \ref{eq:Ineqs} and \ref{eq:Sharelimits} lead to specific capacity factors for flexible sources, described in \cite{Mercure2011TWP}.

These constraints must be respected in order to maintain grid stability, and thus must be applied to the shares equation (eq.~\ref{eq:Shares}). It is done with the interpretation that investment into specific technologies slows down when investors perceive a risk of seeing their capacity unused for grid stability issues, and restrict their technology choices. It is thus through the investor fear of stranded assets that the limiting naturally occurs, an issue of current significant importance for renewable systems (see for instance the report by \cite{Poyry2009}). Thus we can define a probability of reducing investment $g(S_i,\hat{S}_i)$ near the limit given by eq.~\ref{eq:Sharelimits}, and associated cumulative probability distribution of investing given whether the system has past that limit or not, $G(S_i,\hat{S}_i)$, which goes to zero past the limit. Systems can have either an upper limit or a lower limit depending on their nature. The central term of the shares equation can thus be rewritten as 
\beq
A_{ij} F_i(\Delta C_{ij})G_i^{max}G_j^{min} - A_{ji} F_j(-\Delta C_{ij})G_j^{max}G_i^{min},\nn
\label{eq:LimitedSharesEquation}
\eeq
where $G^{max}_i$ represents the probability of investing given $S_i$ and an upper limit $\hat{S}_i$ for technology $i$, if it has one (i.e. if $i$ is a variable or base load system), while $G^{min}_i$ represents a lower limit, if applicable (i.e. if $i$ is a flexible system). Near an upper limit, units of shares cannot flow into a category but can flow away from it. Conversely, near a lower limit share units cannot flow away from a category but can flow in. These constraints, when introduced into the shares equation, lead naturally to a splitting of the electricity generation technology market when the system is close to the limit, for instance constrained by the amount of flexibility. In such a split market, submarkets forms in which competition occurs at different cost levels.\endnote{A submarket for peak load electricity almost always exists. A submarket for variable renewables can also emerge out of various subsidy schemes. Thus three submarkets can exist.}

\subsection{Simple numerical examples with four categories \label{sect:examples}}
\begin{figure}[hbt]
	\begin{center}
		\includegraphics[width=1\columnwidth]{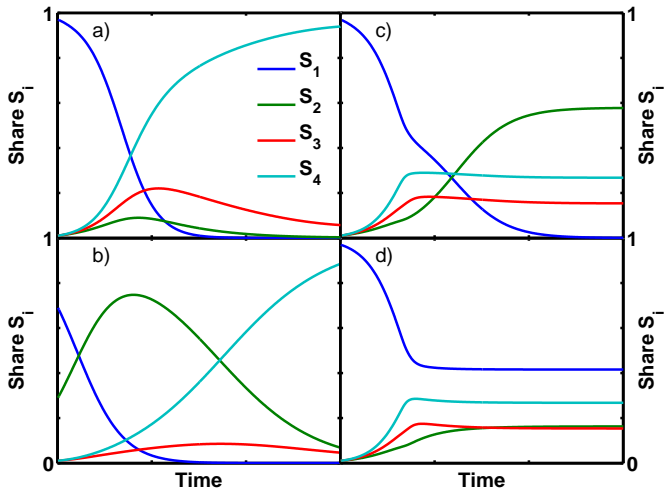}
	\end{center}
	\caption{Examples of calculations using the shares equation in different situations, where constant costs are assumed, with $C_1 = 1.3 C_4$, $C_2 = 1.1 C_4$ and $C_3 = 1.05 C_4$, associated with $S_1$, $S_2$, $S_3$ and $S_4$. $a)$ Starting shares of 97\% for $S_1$, 1\% for the others, no limits are imposed. $b)$ Starting shares of 69\% for $S_1$, 29\% for $S_2$ and 1\% for the others, no limits are imposed. $c)$ Shares as in $a)$ using upper limits of 20\% and 30\% for $S_3$ and $S_4$. $d)$ Shares as in $a)$ using upper limits as in $c)$ with the additional lower limit of 40\% for $S_1$.}
	\label{fig:Figure3}
\end{figure}

We present a brief exploration of market competition as modeled by the shares equation, shown in figure~\ref{fig:Figure3}. The system is defined with four technology categories, the shares of which sum to one, denoted $S_1$ to $S_4$, using identical cost distributions profiles but with different median values, and $A_{ij} = 1$. The constant median cost values are given in terms of the lowest value $C_4$: $C_1 = 1.3 C_4$, $C_2 = 1.1 C_4$ and $C_3 = 1.05 C_4$. Two situations are given without share limits, in $a)$ and $b)$, and two with some limiting, in $c)$ and $d)$. 

In $a)$, most of the market is given to $S_1$, the most expensive technology, with 97\% and 1\% for the other three. We observe that with time, $S_1$ gradually loses the market while the three competitors increase their share exponentially, consistent with the small share limit of a logistic equation. However, competition later arises between these three, where only one, the less expensive $S_4$, wins the whole market. Therefore, with an unconstrained shares equation with constant costs, the technology with lowest cost always eventually wins. This corresponds to an equilibrium state which occurs at a time longer than several times the natural time constant associated with these technology changes, of order of several decades. In $b)$, a larger share is given to $S_2$. We observe that $S_2$ initially increases its share despite the fact that two other technologies have lower costs. This is consistent with a population growth phenomenon, where a larger population, albeit having a lower birth rate, is able to grow faster than a small population. Given enough time, however, the technology with lowest cost still wins in the equilibrium limit. These examples demonstrate the equivalence between a pairwise comparison of technologies performed with a differential equation and a simultaneous comparison of all technologies by exhibiting true multiple interactions, provided by the coupled differential form of the shares set of equations.

We now introduce limiting, first in $c)$, where upper limits of 20\% and 30\% were given to $S_3$ and $S_4$. We observe that these saturate at their limit, and the remaining market is taken by $S_2$, the second most expensive technology. Therefore, with limits, it is not the least expensive that wins the whole market, and we obtain a stable heterogenous system. Adding a lower limit of 40\% for $S_1$, we see that all limited technologies become locked to their limit and that the unlimited technology $S_2$ becomes stable.

\section{Natural resource use and depletion \label{sect:resources}}

\subsection{The cost-supply probability distribution framework}
\begin{figure*}[hbt]
	\begin{center}
		\includegraphics[width=1.4\columnwidth]{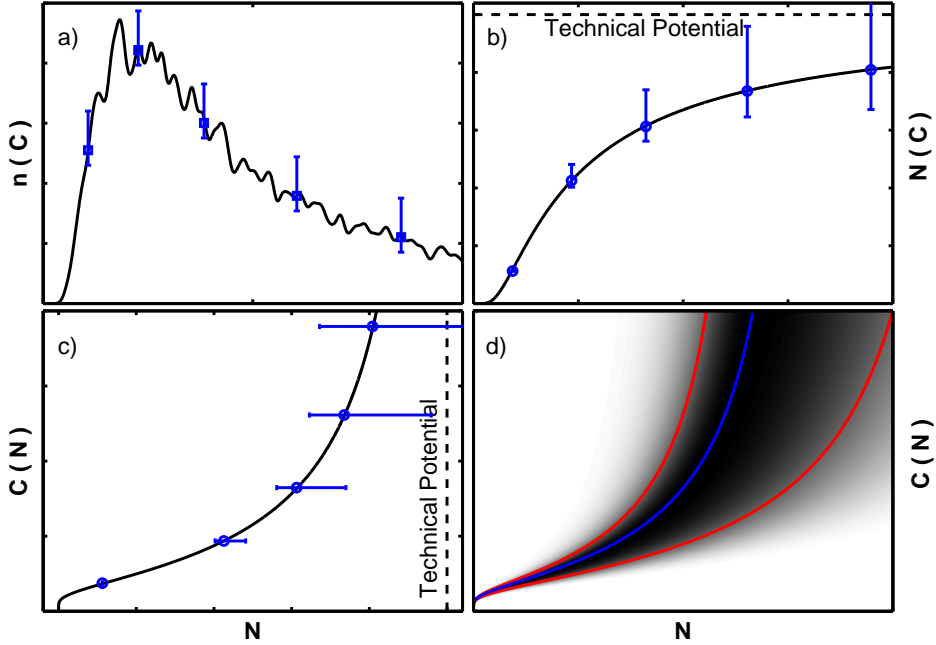}
	\end{center}
	\caption{$a)$ Sketch of a hypothetical density $n(C)$ of cost ranked energy or energy flow units available within cost range $C$ up to $C + \Delta C$, with associated uncertainty indicated using error bars. $b)$ Number of energy or energy flow units $N(C)$ available for a cost of extraction below the arbitrary value $C$. $c)$ The cost $C(N)$ of extracting an additional energy or energy flow unit given that $N$ were already exploited, commonly called the cost-supply curve $d)$ Cost-supply curve defined as a probability distribution, where the red curves indicate the 95\% confidence level range and the blue curve the most probable cost-supply curve.}
	\label{fig:Figure4}
\end{figure*}

Regional energy landscapes are defined by the local availability of natural resources. Local resources can be assessed and ranked in order of cost of extraction. Naturally, a regional electricity sector will develop according to these, unless trading with neighbouring regions provides a more cost-effective solution. The feasibility of a global transition from fossil fuels to sustainable energy sources strongly depends on assessments of global natural resources. These, as given in for instance the World Energy Assessment \citep{WEA2001, WEA2004} are incomplete when they do not provide associated cost distributions. When natural resources occur with low energy densities, in some cases, the underlying assumption to the idea of harnessing such sources implies unrealistically large capital investments, such as, for instance, installing large numbers of wind turbines on many low wind speed sites. The resulting expected aggregated energy production may be very large, but in reality no investor would be willing to undertake such an unprofitable venture that involves large production costs per unit energy.

The framework of cost-supply curves, widely used after the work of \cite{Rogner1997} and taken again by \cite{Hoogwijk2004, Hoogwijk2005, Hoogwijk2009, HoogwijkThesis} for their assessments of wind, solar and biomass global energy potentials and used for instance in the TIMER model \citep{IMAGE} addresses this question. We extend here this approach to include the missing treatment of uncertainty and the availability and evolution of knowledge with respect to natural resources. We use these definitions as one of the central aspects driving the shares equation in FTT:Power, in order to produce meaningful regional energy landscapes. This framework is depicted in figure~\ref{fig:Figure4}. 

We define the density $n(C)$ of units of energy, for non-renewable sources, or units of energy flows, for renewable sources, as a function of their cost of extraction $C$ (see figure~\ref{fig:Figure4} panel $a$). This corresponds to a histogram of the number of energy units that can be extracted for a cost within a certain range $\Delta C$. Data points in this density function possess an uncertainty value (schematically depicted using error bars), corresponding to the lack of exact knowledge over the number of energy units that can be extracted within a particular cost range. The number of units that can be extracted for a cost lower than the arbitrary value $C'$ is the integral of the density $n(C)$ up to $C'$, that we denote $N(C')$ (panel $b$). This function converges towards a unique value at high $C'$, which corresponds to the total technical potential, and is the total area under the density $n(C)$. Note that the associated uncertainty grows cumulatively. This relationship can be inverted in order to express the cost of extracting an additional unit of energy or energy flow after $N$ have already been exploited, denoted $C(N)$ (panel $c$). This cost-supply relationship diverges at the total technical potential. The divergence corresponds to, for example, the installation of a diverging number of wind turbines on a very large number of sites with vanishingly small average wind speeds, but nevertheless giving a finite amount of energy flows. Using uncertainty values, one can in principle determine three cost-supply curves in order to define a confidence range for where the real actual cost-supply curve may lie. This replaces the traditional cost-supply curve by a probability distribution that reflects imperfect knowledge. Panel $d$ depicts this, where the red curves determine the 95\% confidence level region in cost-quantity space, and the blue curve corresponds to the most probable set of values.

Uncertainty in the determination of natural resource availability is notable in the case of fossil fuel reserves and resources. \cite{Rogner1997} paints a clear picture of the process of fossil reserve expansion. In this view, discoveries of deposits of hydrocarbons expand known resources, but their costs of extraction are highly uncertain.\endnote{Note that uncertainty values over cost of extraction can be transformed into uncertainty over amounts of energy within well defined cost ranges; it only requires a redistribution of the same data.} However, with technological learning and additional exploration, knowledge over cost values improves but costs also gradually decrease, resulting in a flow from uncertain and costly resources to known and economic reserves. Thus, as time progresses into the future, the cost-supply curves of fossil resources become better and better defined, an effect that can strictly only be described by a probability distribution defined now based on current knowledge, which determines the set of probable future paths.

\subsection{Cost-supply curves in FTT:Power}

Cost-supply curves are necessary to define the level of use of each natural resource in a model for electricity production such as FTT:Power. The underlying assumption is that within a specific world region, units of a particular resource will be used in order of cost to the best of current knowledge. In the case of wind power for instance, wind sites with particularly good wind speed distributions and low turbulence properties are likely to be developed first, furthermore in areas where they are likely to be socially accepted and where land is affordable. All these aspects should therefore be included in such assessments, an effort which has been done by Hoogwijk $et$ $al.$ in an extensive assessments of global wind potentials \citep{Hoogwijk2004}. However, the level up to which resources are used depends entirely on the availability and cost of all alternatives which can be used to produce electricity in that region in order to meet the demand. 

This is precisely what is done in FTT:Power, where cost supply curves are used for every calculation of the LCOE as defined by eq.~\ref{eq:LCOE}. This is an aspect of central importance, since the complete set of cost-supply curves maintain the calculation within boundaries defined according to a complete set of total technical potentials and thus yields reasonable results of natural resource use. Additionally, a marginal cost of electricity production can simply be derived using the LCOE averaged over all technologies, $\sum_i S_i \times LCOE_i$ ($S_i$ being the shares of capacity), and used with additional information regarding taxes and profit margins to determine the price of electricity. The price of energy carriers being important for the global economy, as represented for instance in E3MG, this aspect of the model has a deep influence over future projections of energy demand and supply and associated effects on economic activity. Moreover, cost-supply curves inherently track closely resource depletion.

The set of cost-supply curves in a specific world region define a marginal cost of electricity production, and therefore results in all world regions having different prices of electricity. In a real world, this inevitably results in electricity trading between regions which have their grids connected, from resource poor to resource rich regions. Examples of this are the sale of important amounts of Canadian or Scandinavian hydroelectricity in American or European markets respectively. Thus except for closed energy markets, the level of use of natural resources for electricity production determined through cost-supply curves of world regions should therefore additionally depend on each other. The details and effects of energy trading between the 20 E3MG regions within FTT:Power are complex and will be explored elsewhere. Moreover, with the introduction of other sectoral models of technology of the FTT family into E3MG, competition for resources will arise, an aspect that will be discussed in later work.

FTT:Power uses deterministic cost-supply curves, or in other words, a single curve amongst all possible paths defined by the probability distribution of figure~\ref{fig:Figure4} $d)$. Although the decision making underlying the shares equation (eq.~\ref{eq:Shares}) should include an assessment of risk and uncertainty in investment, a subject which is not the scope of the present work but will be the theme of a subsequent paper, the uncertainty over natural resources does not influence the choices of investors directly. This is due to the fact that the cost-supply curve always has a vanishing uncertainty at the current level of natural resource use, exemplified by a well defined price. It is the level of future use which is uncertain, generally beyond the typical time horizon for capacity investment. Thus, at every time step of the simulation, resource availability and costs as seen by investors should be well defined quantities, without uncertainty. It is the number of possible futures which is large. The determination of uncertainty over predictions therefore requires runs of the model using large numbers of possible deterministic cost-supply curves derived from probability distributions akin to that in figure~\ref{fig:Figure4} $d)$, for all sectors. Monte-Carlo techniques are thus the appropriate way forward in order to define confidence levels attributed to future predictions, an aspect we treat as essential in the formulation of this model, and will be explored in detail elsewhere.

\section{Results \label{sect:results}}
\subsection{Model parameters and assumptions}

\begin{figure*}[t]
	\begin{center}
		\includegraphics[angle = -90, width=1.7\columnwidth]{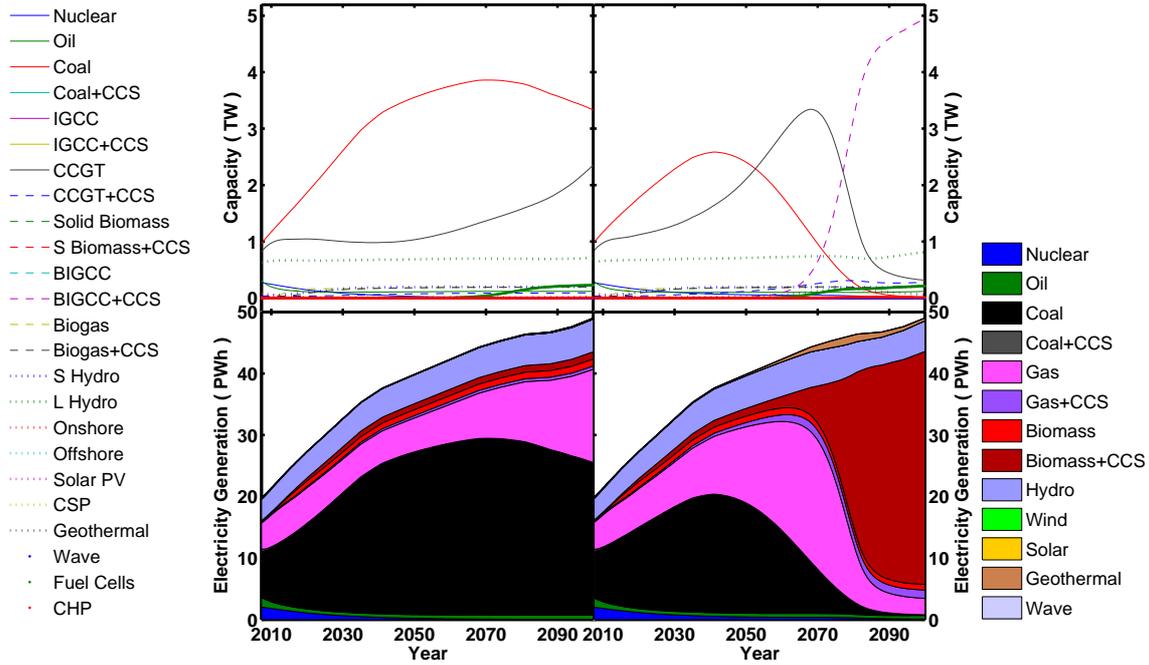}
	\end{center}
	\caption{Model results for two sets of assumptions, a baseline (left panels) and a mitigation (right panels) scenarios. The top panels show the capacities $U_i$ for 24 energy technologies, while the bottom two show electricity generation $G_i$ as areas that sum up to the total electricity demand. For $G_i$, related categories were aggregated into a smaller set for clarity. }
	\label{fig:Figure5}
\end{figure*}

The implementation of FTT:Power which was used to produce the results presented here was made with the definition of 24 energy technologies and one single world region, along with global cost-supply curves for 13 types of resources, some of which are used by several competing energy technologies. A table summarising these is provided in the appendix. As a standalone model, it essentially requires a curve for energy demand and starting values for the capacities $U_i$, and calculates energy capacity and generation from eq. \ref{eq:G} and \ref{eq:U}, with eq. \ref{eq:Shares} driving the changes in the energy mix. For the sake of exploring the properties of FTT:Power and presenting results which stem from the simplest assumptions, we have excluded any feedback with the global economy, from electricity demand and through electricity prices, by running the model by itself without E3MG. Therefore, we used an exogenous electricity demand curve. 

Components of the LCOE, given in the appendix, were derived from both recent IEA data \citep{IEAProjCosts} and from a set of global cost-supply curves \citep[a detailed analysis is forthcoming, in ][]{Mercure2012}. These curves were constructed from an ensemble of sources in order to produce a complete set \citep{IAEA2009, Rogner1997,Hoogwijk2004,Hoogwijk2005,Hoogwijk2009,HoogwijkThesis,Themelis2007,Mock1997, IJHD2011, Lako2003}. In the case of Biomass, the global potential depends highly onto assumptions regarding world population and global food demand, and the results of \cite{Hoogwijk2009} for the B1 SRES scenario were used. For this paper, as a single region global model, the assumption is made of perfect electricity trading worldwide, such that every type of resource is available everywhere, an assumption which is done here for the convenience of using a single world region for the purpose of the demonstration of the validity of the model equations, but which will be dropped after the subsequent definition of cost-supply curves for 20 world regions that coincide with those of E3MG, an aspect of FTT:Power which is not the focus of the present work . The results that we present, although highly aggregated, are of interest in themselves and in how they depict the technology substitution processes through which a transition towards a low carbon electricity sector may occur in the future.

We present here results for two simple scenarios, a baseline where no mitigation effort is made but CO$_2$ pricing exists, and a scenario where emissions are decreased through reducing the number of allowances such that the price for allowances increases exponentially (in real 2008 dollars). We thus assumed a constant price of carbon of 22\$/t CO$_2$ in the baseline scenario, while in the mitigation scenario the price of carbon is the same initially but increases by 1\% per year. This was done in order to have an identical starting point, such that a diverging evolution of the nature of new energy capacity can be observed, entirely attributable to the difference in the price of carbon.\endnote{Without 22\$/t CO$_2$ in the baseline, coal power plants dominate entirely the energy mix, as one would expect.} With FTT:Power as part E3MG, the price of carbon will be endogenous. The price of carbon is fed back into the system through the carbon component of the LCOE (eq.~\ref{eq:LCOE}), where a high discount rate of 10\% was used.\endnote{The choice of the discount rate influences highly technology choices by setting the strength of the preference for delaying costs to the future. Thus, with a high discount rate, technologies with lower initial investment costs relative to other components of the LCOE are preferred. Effects of the discount rate are discussed in \cite{Mercure2011TWP}. All discount rates can in principle be set separately exogenously in this model.} Therefore, an economic incentive to depart from heavily emitting technologies and edge towards more efficient systems is created. For simplicity, however, no feed-in tariffs or subsidies were assumed for any of the renewable categories, and therefore wind and solar systems do not appear in the results. Both scenarios use an identical demand function, which was derived from that of the New Policies Scenario of the World Energy Outlook 2010 \citep{IEAWEO2010}. It was extended beyond 2035 assuming a slightly decreasing growth rate, chosen as a good example rather than for accuracy. In the real world, the demand is expected to be reduced in mitigation scenarios where increases in production costs due to an increasing price of carbon or technology switching is likely to be passed on to the consumer through increasing electricity prices. Such an effect would constitute a direct result of feedback between FTT:Power and E3MG.

\subsection{Model Outputs}

Figure~\ref{fig:Figure5} presents the results for the capacity $U_i(t)$ (top panels) and the power generation $G_i(t)$ (bottom panels), the latter being represented as coloured areas which sum up to the total electricity demand. The baseline is shown in the left panels, while the mitigation scenario is shown in the right panels. The baseline scenario is dominated by coal fueled power stations, except for a sizable amount of hydroelectricity, which is seen to converge towards a capacity of about 0.7~TW worldwide. This is due to the fact that a large fraction of the world's total hydroelectric potential is already developed, and that costs rise very rapidly with additional dams being built. Thus, this particular single world region calculation predicts a scarcity of suitable economic sites for additional hydroelectric projects.\endnote{Note that here regional interactions would be important, since in areas where alternatives are also expensive, expensive hydroelectric dams could still be built. This will be addressed in a subsequent paper.} The remaining demand is covered mostly by coal power stations and combined cycle gas turbines (CCGT).

The mitigation scenario exhibits more richness with a larger variety of substitution events triggered by the increasing price of carbon. Hydroelectricity appears hardly affected by the difference in assumptions. However, as the price of carbon increases, coal is gradually phased out, peaking in around 2040, where CCGT takes the lead. This is due to the much increased efficiency of CCGT, which emits less GHGs per unit of energy compared to conventional coal. Simple CCGT systems are however in turn gradually phased out from 2070 onwards to be replaced by  systems with negative emissions, biomass integrated gasification combined cycle (BIGCC) power plants with carbon capture and storage (CCS) technology. Since such systems capture CO$_2$ from the atmosphere through growing biomass crops, and subsequently sequestrate it into permanent storage, they effectively emit a negative amount of CO$_2$ into the atmosphere. Such a process should in principle be associated to a generation of emissions allowances which can be sold onto the carbon market, and thus generate income, resulting in a negative cost component in the LCOE. 

CO$_2$ emissions reductions resulting from this mitigation scenario are substantial. In comparison to the baseline, power sector emissions peak at a value of around 20~GtCO$_2$/y in 2045, compared to a stabilisation at 23~GtCO$_2$/y for the baseline. In 2080, complete decarbonisation of the power sector is achieved, as opposed to emissions of 27~GtCO$_2$/y for the baseline. After this date, emissions rapidly become negative through CO$_2$ sequestration of emissions generated from the combustion of biomass fuels. Thus, cumulative emissions reach a peak in 2080 of 1000~GtCO$_2$ relative to 2010, and decrease afterwards. We emphasise our use of an energy demand curve derived from the New Policies Scenario of the IEA \citep{IEAWEO2010}, which does not include significant demand reduction measures in the mitigation effort. Emissions are expected to be reduced further by such measures through a decreased global supply of electricity. However, large uncertainties arise associated with the electrification of various other energy services which do not rely on electricity, such as heat and transport, that also have the potential to alter significantly the global demand for power. The electrification of such services could lead to major emissions reduction if the power sector is itself already decarbonised, but could also lead to major increases if it is not, through large efficiency losses (for instance by replacing auto fuels by coal based electricity).

\subsection{The energy technology ladder}

\begin{figure*}[t]
	\begin{center}
		\includegraphics[width=1.3\columnwidth]{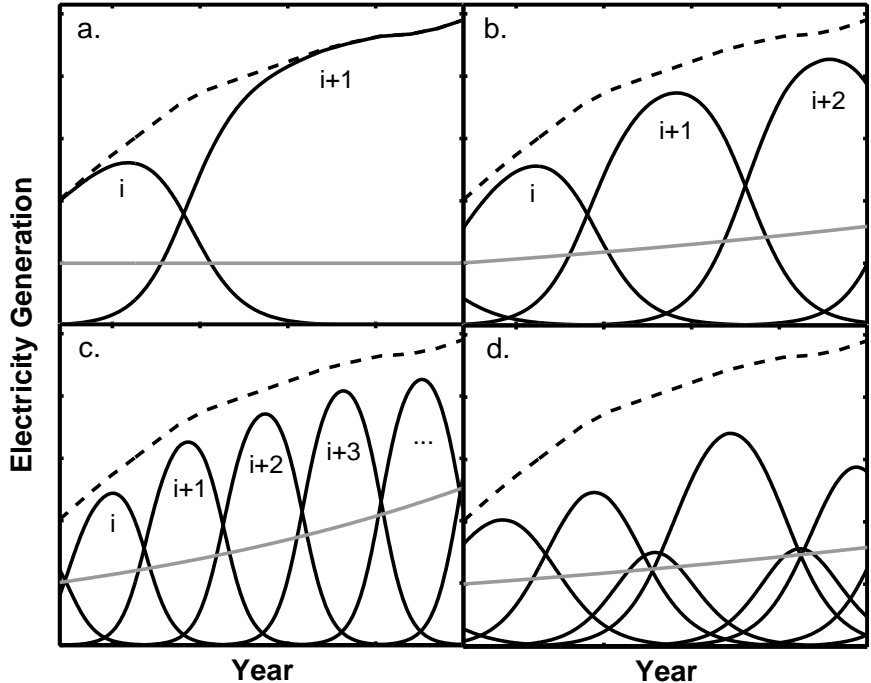}
	\end{center}
	\caption{Sketch of the concept of the energy technology ladder, where energy technologies are gradually replaced as the price of emissions allowances increases. The electricity generation by technology is shown as solid lines, the total electricity demand is the dashed line, while the price of emissions allowances is shown as a grey curve. Technology $i$ is replaced by $i+1$, which in turn is replaced by $i+2$ and so on, where at each point in time the mix of technology possesses a marginal cost of abatement, the cost of an additional substitution towards technologies with lower emissions, equal to the price of emissions allowances. $a.$ Without carbon pricing, technological learning may still generate a slow technological transition. $b.$	Low rate of increase of the price of allowances. $c.$ High rate of increase. In this case, many more intermediate technologies are used. $d.$ More realistic picture with emission factors not evenly distributed and varying substitution rates.}
	\label{fig:Figure6}
\end{figure*}

The technological transitions presented above in the mitigation scenario are typical of the type of effect that are observed in FTT:Power when technological change is favoured by policy assumptions such as an increasing price of carbon. These results are strikingly similar to technological transitions observed historically as reported by Marchetti and Nakicenovic \citep{Marchetti1978} and by Grubler $et$ $al.$ \citep{Grubler1999}. Several such technological transitions should be expected in a high mitigation scenario aimed for instance at the decarbonisation of the energy sector. Due to the gradual nature of learning and the continuous increase of the price of carbon, this transformation is likely to be eased by the appearance of transient states. Rather than jumping directly to the lowest emission technologies immediately, the system adopts intermediate solutions which it subsequently phases out. Due to the lifetime of systems, it cannot perform the transitions faster than the time it takes for power plants to come to the end of their useful life. Power stations are gradually replaced by other types with lower emissions, which are in turn also replaced by yet more efficient ones, and so on. We name this effect the energy technology ladder. 

We sketch this phenomenon graphically in figure~\ref{fig:Figure6}. We first rank energy technologies in terms of their factor of GHG emissions per unit energy and denote them $i$, $i+1$, $i+2$, and so on, where $i$ corresponds to currently dominating heavy emitting coal power stations. In all panels, the dashed line represents the total electricity demand, while the solid lines are the power generation by individual technologies. The grey line represents the price of carbon in a different set of arbitrary units. In panel $a$, without an increase in carbon price, one slow transition may possibly take place, driven by learning and, for instance, conversion efficiency and relative price of fuels. In $b$, however, the slowly increasing price of carbon triggers several gradual substitutions in the power sector similarly to those observed in the mitigation scenario shown above, where two substitution processes are observed. In $c$, a faster increase of the price of carbon allows the appearance a larger number of systems, each with ever lower emission factors, going through four substitution processes, and results in a final system with lower total emissions. In $b$ and $c$, equally spaced emission factors and costs lead to evenly distributed substitution events as a function of time. Emission factors are in reality not evenly distributed, and rates of substitution are not expected to be identical due to different system lifetimes and lead times. Real power systems are likely to follow the more realistic picture shown in $d$, where more than single substitution events may occur simultaneously and at different rates.

The energy technology ladder is an effect which emerges from the equations underlying FTT:Power in a variety of sets of assumptions, whenever technological change is favoured. It is a general result that stems from the combination of technological learning as given by experience curves (eq.~\ref{eq:ExpCurve}) and a shares equation based on a logistic set of differential equations (eq.~\ref{eq:Shares}), and leads to classic sigmoid ($S$-shaped) technological transitions. We stress that this property is not an equilibrium property, and that not all technological substitutions made economic by the price of carbon occur at any one time. This is due to the dynamic nature of the shares equation, which takes fundamental account of sector growth time constants. Therefore, by including the dynamics of growth, one cannot assume an equality to exist between the price of carbon and the marginal cost of abatement, since an equilibrium is never reached. 

We moreover consider that other sectors of the economy could be modeled in a similar way, and thus see logistic technological transitions and technology ladders. One of the sectors of interest is transport, in which the world may see similar types of technological transitions between, for instance, petrol based vehicles and other types such electric cars or systems running on biofuels. As shown by \cite{Grubler1999} for transport networks and infrastructure, this sector is likely to follow the same logistic behaviour with its own set of time constants. Technology transitions in the transport sector are likely to have massive influence on the energy sector as a whole, and should thus be included in any energy model in order to be complete, since both systems are highly correlated to each other. Such a model for transport is under consideration as a member of the FTT family of bottom-up models connected with E3MG, FTT:Transport.

\section{Conclusion}

We have introduced in this work the model FTT:Power for the simulation of global power systems based on a set of logistic differential equations and induced technological change. The model was designed to be integrated as one member of the FTT family of bottom-up models into the global macroeconometric model E3MG. Induced technological change arises with learning, which in combination with logistic differential equations, leads to classic irreversible and path dependent $S$-shaped technological transitions. Logistic differential equations offer an appropriate treatment of the times involved in the transformation of technology for the power sector; in particular it leads to exponential rises at low penetration, and saturation at high penetration. Competition between technologies occurs through a pairwise comparison of the LCOE, which we have demonstrated is equivalent to a simultaneous comparison of the LCOE of all options. Constraints related to the electrical grid and its properties are expressed as limits on shares of technologies. Restrictions on natural resource cost-effective use and depletion are provided by the use of probability distributions associated with cost-supply curves, providing an appropriate treatment of uncertainty in resource assessments. 

An example of a calculation performed with FTT:Power for two simple global scenarios, a baseline and a mitigation context, using global parameters and cost-supply curves, is given in order to explore the properties of the model. The baseline features constant carbon pricing while the mitigation scenario sees a carbon price increase by 1\% per year. Baseline assumptions result in a resource limited amount of hydroelectricity, the remaining demand being mostly covered by coal-powered electricity generation. The mitigation scenario yields similar results for the first decade, but gradually diverges towards lower GHG emissions and two major technological transitions. The first occurs between coal and gas, while the second transforms gas into biomass based electricity generation while sequestrating CO$_2$, resulting in a negative component in the LCOE due to negative emissions. The observed succesive transformations of the energy sector are akin to classic observed logistic technology transitions, and are driven by the pricing of GHG emissions. We name this effect the technology ladder. Higher rates of increase of the price of carbon lead to faster and richer transformations and to lower GHG emissions. 

\section*{Acknowledgements}
The author would like to acknowledge T. S. Barker for guidance and support, P. Salas, for help with data and model testing, A. Anger (Cambridge), H. Pollitt and P. Summerton (Cambridge Econometrics) for informative discussions. I would also like to thank A. Gr\"ubler at IIASA for generating lively debates, prompting a thorough verification of the results. This work was supported by the Three Guineas Trust.

\appendix
\section{Description of technologies used in the model}

\begin{table}[h]\footnotesize
\begin{center}
		\begin{tabular*}{1\columnwidth}{@{\extracolsep{\fill}} l|r r r r r r r }
			\hline
			\multicolumn{7}{ l }{Model technologies and parameters}\\
			\hline
			Technology & \multicolumn{4}{ l }{Cost components}  &\multicolumn{3}{ r }{\cite{IEAProjCosts}}\\
			Name & $I$ & $F$ & $O\&M$ & $\alpha$ & $\tau$ & $t$ &$b$\\
			& \scriptsize{\$/kW} & \scriptsize{\$/MWh} & \scriptsize{\$/MWh} & \scriptsize{t/GWh} & y & y \\
			\hline
			\hline
			Nuclear 		&3739	&9.33	&14.23	&0		&60	&7	&0.086\\
			Oil  			&1139	&207	&20.53	&586	&40	&4	&0.014\\
			Coal 		&2134	&20.01	&6.38	&852	&40	&4	&0.044\\
			Coal + CCS 	&3919  	&20.81	&13.93	&97		&40	&4	&0.074\\
			IGCC 		&3552	&18.60	&9.36	&852	&40	&4	&0.044\\
			IGCC + CCS 	&4194	&18.52	&11.94	&97		&40	&4	&0.074\\
			CCGT	 	&1047	&60.08	&4.51	&354	&30	&2	&0.059\\
			CCGT + CCS 	&2269	&66.05	&5.96	&43		&30	&2	&0.074\\
			Solid Biomass 	&4491	&44.10	&10.09	&0		&40	&4	&0.074\\
			S. Biom. + CCS	&6277	&44.10	&10.09	&-981	&40	&4	&0.105\\
			BIGCC	 	&3552	&44.10	&9.36	&0		&40	&4	&0.074\\
			BIGCC + CCS 	&4194	&44.10	&11.94	&-981	&40	&4	&0.105\\
			Biogas	 	&2604	&26.50	&24.84	&0		&30	&2	&0.074\\
			Biogas + CCS 	&3826	&26.50	&24.84	&-376	&30	&2	&0.105\\
			Tidal	 		&2611	&0		&44.00	&0		&80	&7	&0.020\\
			Hydro		&2138	&0		&5.11	&0		&80	&7	&0.020\\
			Onshore		&1963	&0		&21.26	&0		&25	&1	&0.105\\
			Offshore		&4453	&0		&39.40	&0		&25	&1	&0.136\\
			Solar PV		&5153	&0		&23.73	&0		&25	&1	&0.269\\
			CSP			&5141	&0		&27.59	&0		&25	&1	&0.152\\
			Geothermal	&5286	&0		&18.21	&0		&40	&4	&0.074\\
			Wave		&4770	&0		&51.87	&0		&20	&1	&0.218\\
			Fuel Cells		&5459	&54.46	&49.81	&69		&20	&2	&0.234\\
			CHP			&1529	&55.84	&9.20	&69		&40	&2	&0.044\\
			\hline
		\end{tabular*}
	\caption{List of technologies used in the model with values assumed for its parameterisation. }
	\label{tab:values}
\end{center}
\end{table}

Table \ref{tab:values} provides a complete description of technologies used in the model. Parameters were obtained mainly from statistics performed using values from \cite{IEAProjCosts}. Additional cost uncertainty values used are not shown. $CCS$ stands for Carbon Capture and Storage, Coal stands for all types of coal power stations excluding $IGCC$, $IGCC$ stands for Integrated Gasification Combined Cycle used with coal fuel, $CCGT$ for Combined Cycle Gas Turbine, $BIGCC$ for Biomass IGCC, $PV$ for photovoltaic, $CSP$ for Concentrated Solar Power and $CHP$ for Combined Heat and Power where the heat generated by CCGT power stations is used for industrial purposes, improving the effective efficiency by replacing an equivalent amount of electricity. $I$ refers to investment costs, $F$ to fuel costs, $O\&M$ to Operation and Maintenance, $\alpha$ is the technology specific emissions factor obtained from \cite{IPCCGuidelines}, $\tau$ is the technology lifetime and $t$ is a parameter proportional to the comparative time scale of industry production capacity expansion. $b$ is absolute value of the learning exponent. Proxies or estimates were used for missing values. Starting share, capacity and generation values, as well as the electricity demand were obtained from \cite{IEAWEO2010}. 

\theendnotes





\bibliographystyle{elsarticle-harv}
\bibliography{CamRefs}







\end{document}